\documentclass[a4paper,twoside,11pt]{amsart}
\usepackage{hyperref}
\usepackage{color}
\usepackage{amsmath}
\usepackage{graphicx}
\usepackage[utf8]{inputenc}
\usepackage{amsfonts}
\usepackage{geometry}
\usepackage{amssymb}
\usepackage{latexsym}
\usepackage[english]{babel}
\usepackage{amssymb}               % para cargar las fuentes de la A.M.S.
\usepackage{amsmath}               % para usar letras griegas en negrita
\usepackage{amsthm}  

\newtheorem{theorem}{Theorem}[section]

\newtheorem*{remark}{Remark}
\newtheorem{lemma}[theorem]{Lemma}
\newtheorem{corollary}[theorem]{Corollary}

\newtheorem{proposition}[theorem]{Proposition}

\makeindex
\theoremstyle{plain} % Entornos tipo teorema

\theoremstyle{definition}
\title[Obstructions to the categorification of near-group rings]{Galois-modular obstructions to the categorification of near-group fusion rings}
\author{Víctor Blasco}

\address{Department of Mathematics,
Louisiana State University,
Baton Rouge, LA 70803, USA}

\email{vblasc1@lsu.edu}
\begin{document}
\begin{abstract}
    In this article we establish new obstructions to the categorification of near-group fusion rings by combining the values of irreducible characters of the Grothendieck ring with the Galois symmetry of the induced unit object $I(1)$ in the Drinfeld center. As a result, we show that a near-group ring of type $G+m$, where $G$ is an abelian group and both $|G|$ and $m$ are odd, can only categorify if  $m=|G|$. 
\end{abstract}
\maketitle
\section{Introduction}
The categorification problem for fusion rings asks which fusion rings can be realized as the Grothendieck rings of fusion categories. A particularly interesting instance of this problem concerns near-group fusion rings, whose basis consists of a finite group $G$ of invertible elements together with a single non-invertible element $X$. More concretely, a near-group fusion ring of type $G+m$  has \(G\cup\{X\}\) as a basis, with fusion rules
\[
gh=gh,\qquad gX=Xg=X,
\]

\[
X^2=\sum_{g\in G}g+mX,
\]
\medskip

Despite the simplicity of their fusion rules, determining necessary and sufficient conditions for near-group rings to admit a categorification remains a difficult problem; and a complete classification is still unknown, but some partial results have been achieved. The case $m=0$ was completely classified in \cite{TY}. For $G$ abelian and positive $m$, \cite{Si} established the general restriction $m\geq |G|-1$ and proved that, when $m=|G|-1$, categorification is only possible precisely when $G$ is the multiplicative group of a finite field. From an arithmetic perspective, Ostrik \cite{O1,O2} developed obstructions involving formal codegrees and the Drinfeld center. More recently, \cite{Sch} established general numerical restrictions on categorifiability, completely classifying the categorifiable near-group rings associated with elementary abelian $2$-groups, and obtained additional obstructions for cyclic groups of prime order $p\equiv 3\pmod 4$.
\medskip

In this paper, we contribute to this classification problem by proving the following necessary condition: if $G$ is an abelian group of odd order and $m$ is also odd, then $G+m$ can only categorify if $m=|G|$. This generalizes \cite[Proposition 4.4]{Sch}. Our approach builds on the arithmetic and representation-theoretic methods developed in \cite{O2}. More specifically, we study the irreducible characters of the Grothendieck ring and their relationship with the simple constituents of the induced unit object in the Drinfeld center. By exploiting the Galois symmetry of these constituents and using the FS-indicators defined in \cite{NS}, we obtain additional restrictions on the possible categorifications, leading to our main non-existence result. The argument provides a relatively direct obstruction without requiring heavy computation.
\medskip

The paper is organized as follows. We begin in Section 2 by recalling some basic facts about fusion categories that will be useful later on. In Section 3 we study the eigenvalues of the Casimir operator at each simple object $S$, and their relations with the Casimir element and the characters of the Grothendieck ring of the fusion category. In Section 4 we relate those character values with the $S$-matrix of the Drinfield center. This helps to prove in Section 5 that $I(1)$ is Galois closed, where $I$ denotes the induction functor from $\mathcal{C}$ to $\mathcal{Z}(\mathcal{C})$. Finally, in Section $6$ we prove the result concerning the obstruction for categorification of certain near-group categories by looking at how a particular Galois automorphism permutes the unit object $\mathbf{1}$ in $\mathcal{Z}(\mathcal{C})$.
\bigskip

\section{Preliminaries}
Let $\mathcal{C}$ be a fusion category over $\mathbb{C}$, and let $\operatorname{Irr}(\mathcal{C})$ denote the set of isomorphism classes of its simple objects. The fusion rules of $\mathcal{C}$ are encoded in its Grothendieck ring $K(\mathcal{C})$, whose multiplication is induced by the tensor product. The unit object is denoted by $1$. The Frobenius--Perron dimension is a ring homomorphism $\operatorname{FPdim}:K(\mathcal{C})\to\mathbb{R}$, positive on simple objects, and $\operatorname{FPdim}(\mathcal{C})
=\sum_{X\in\operatorname{Irr}(\mathcal{C})}\operatorname{FPdim}(X)^2$.
\medskip

A fusion category $\mathcal{C}$ is called pivotal if it is equipped with a monoidal natural isomorphism $j:\operatorname{Id}_{\mathcal{C}}\to(-)^{**}$. A pivotal structure is said to be spherical if the left and right pivotal traces coincide for every endomorphism $f:X\to X$, that is, $\operatorname{Tr}_L(f)=\operatorname{Tr}_R(f)$. We will assume our categories to be always spherical. In this case, the categorical dimension of an object $X$ is $\dim(X):=\operatorname{Tr}_L(\operatorname{id}_X)
=\operatorname{Tr}_R(\operatorname{id}_X)$, and the global dimension of $\mathcal{C}$ is given by $
\dim(\mathcal{C})=\sum_{X\in\operatorname{Irr}(\mathcal{C})}\dim(X)^2$. The category $\mathcal{C}$ is called pseudo-unitary if $\dim(\mathcal{C})=\operatorname{FPdim}(\mathcal{C})$. \bigskip

%The complexified Grothendieck ring is a finite-dimensional semisimple algebra. Associated with each of its irreducible representations $E$, with character $\chi_E$, is a positive real algebraic number $f_E$, called the *formal codegree* of $E$. These invariants are closely related to the structure of the Drinfeld center $\mathcal{Z}(\mathcal{C})$, whose objects are pairs consisting of an object of $\mathcal{C}$ and a compatible half-braiding. The forgetful functor $F:\mathcal{Z}(\mathcal{C})\to\mathcal{C}$ admits a two-sided adjoint $I:\mathcal{C}\to\mathcal{Z}(\mathcal{C})$, called the *induction functor*, satisfying
The Drinfeld center $\mathcal{Z}(\mathcal{C})$ is the braided fusion category whose objects are pairs $(X,c_{X,-})$, where $X
$ is any object in $\mathcal{C}$ and $c_{X,-}$ is a half-braiding on $X$, namely, a natural family of isomorphisms
$$
c_{X,Y}: X\otimes Y\longrightarrow Y\otimes X
$$
satisfying the hexagon identity.
\bigskip

The forgetful functor $F:\mathcal{Z}(\mathcal{C})\to\mathcal{C}$ sends each pair $(X,c_{X,-})$ to its underlying object $X$. It admits a two-sided adjoint $I:\mathcal{C}\to\mathcal{Z}(\mathcal{C})$, called the induction functor. Moreover, for every $X\in\mathcal{C}$, we have

$$
F(I(X))\cong
\bigoplus_{Y\in\operatorname{Irr}(\mathcal{C})}
Y\otimes X\otimes Y^*.
$$

For spherical $\mathcal{C}$, its Drinfeld center is a modular tensor category, equipped with a ribbon twist $\theta$. On a simple object $A$, the twist acts by a root of unity $\theta_A$. The modular data of $\mathcal{Z}(\mathcal{C})$ is described by its (unnormalized) $S$- and $T$-matrices, given by

$$
S_{A,B}=\operatorname{Tr}(c_{B,A}\circ c_{A,B}),
\qquad
T_{A,B}=\delta_{A,B}\theta_A
$$

%where the $S$-matrix is unnormalized.

The ribbon twist is also related to the Frobenius--Schur indicators $\nu_n(X)$, defined as traces of canonical operators on $\operatorname{Hom}(1,X^{\otimes n})$. They satisfy

$$
\nu_n(X)=\frac{1}{\dim(\mathcal{C})}
\sum_{A\in\operatorname{Irr}(\mathcal{Z}(\mathcal{C}))}
[F(A):X]\,\dim(A)\,\theta_A^n.
$$
\medskip

when the category is spherical. See \cite{EGNO}, \cite{NS} for further details.

\section{Character values and determinant}

Let $K(\mathcal{C})$ be the Grothendieck ring of a spherical fusion category $\mathcal{C}$. We will identify objects in $\mathcal{C}$ with their classes in $K(\mathcal{C})$ in general. Write $R=K(\mathcal{C})\otimes_{\mathbb{Z}}\mathbb{C}$.  For each irreducible representation $(E,\rho_E)$ of $R$, let $\chi_E$ denote its character (i.e. $\chi_E=\operatorname{Tr}\circ\rho_E$). Following \cite{O1}, we consider the central element

$$
\alpha_E=\sum_{X\in\operatorname{Irr}(\mathcal{C})}\chi_E(X)X^*\in R.
$$

The element $\alpha_E$ acts by a scalar $f_E$ on $E$ and by zero on every other irreducible representation. This scalar is called the formal codegree associated with $E$. In particular,

$$
f_E=\frac{1}{\dim(E)}
\sum_{X\in\operatorname{Irr}(\mathcal{C})}
\chi_E(X)\chi_E(X^*).
$$

Formal codegrees are positive real algebraic integers greater or equal than $1$. In particular, both FPdim$(\mathcal{C})$ and dim$(\mathcal{C})$ are codegrees (they are the same codegree when the category is pseudounitary). Moreover, they are closely related to the decomposition of the induced unit object in the Drinfeld center:

$$
I({1})\cong\bigoplus_E \dim(E)\,A_E,
\qquad
\dim(A_E)=\frac{\dim(\mathcal{C})}{f_E},
$$

where the $A_E$ are pairwise nonisomorphic simple objects of $\mathcal{Z}(\mathcal{C})$.
\bigskip

The Casimir operator is the linear map $c:R\to Z(R)$ defined by

$$
c(a)=\sum_{X\in\operatorname{Irr}(\mathcal{C})}XaX^*.
$$

In particular, the central element $c(1)=\sum_{X\in\operatorname{Irr}(\mathcal{C})}XX^*$ is called the Casimir element of $R$. Its action on the irreducible representations of $R$ is closely related to the formal codegrees. More precisely, for every irreducible representation $E$ with character $\chi_E$ and formal codegree $f_E$, we have $
c(x)|_E=f_E\chi_E(x)\operatorname{id}_E$; from where it follows that $c(1)|_E=f_E\dim(E)\operatorname{id}_E$.
\medskip

Thus, the eigenvalues of the action by left multiplication of the Casimir element on $R$ are precisely the numbers $f_E\dim(E)$, with $E$ ranging over the irreducible representations of $R$. The determinant $det_1=\prod_{E} (f_E dim(E))^{dim(E)^2}$ has been considered in \cite{WLL}.

\begin{lemma}
   Let $\mathcal{C}$ be a spherical fusion category with commutative Grothendieck ring and conductor $N=$FSexp$(\mathcal{C})$. Then $det_1$ and $N$ have the same prime divisors.
 \end{lemma}
\begin{proof}
    It follows by \cite{O2} that every formal codegree belongs to the dimension field of the category, and that each $\frac{\dim(\mathcal{C})}{f_E}$ is an algebraic integer; so each codegree divides $\dim(\mathcal{C})$ in $Z[\zeta_N]$. Then, since $\dim(\mathcal{C})$ is always one of the codegrees, and $R$ is commutative, it follows that $det_1=\prod_{E} f_E$. Hence, both $(\dim(\mathcal{C}))$ and $(det_1)$ have the same prime ideal factors in $Z[\zeta_N]$. But Cauchy's theorem for spherical fusion categories (see \cite{BNRW}) tells us that $(\dim(\mathcal{C}))$ and $(N)$ also have the same prime factors. Hence, since both $N$ and $det_1$ are integers it follows that they have the same prime divisors.
\end{proof}
\begin{remark}
    We notice that the condition of commutativity is essential, as the example 4.7 in \cite{Scn} shows. In that case, the determinant is even because the (non-commutative) Grothendieck ring contains a block of size two since the rank is $8$; but the norm of the dimension is $21125$, which is odd.
\end{remark}
\medskip

Let $S$ be any simple object of $\mathcal{C}$, and let $E_1,...,E_r$ be the simple representations of $R$. Call by $\rho_i$ the corresponding ring homomorphism from $R$ to  $End_{\mathbb{C}}(E_i)$, and by $f_1,...,f_r$ the formal codegrees. Now, $c(S)$ is a central element so it acts as an scalar, call it $f_{S,i}$, on $E_i$, meaning that $\rho_i (c(S))=f_{S,i} id_{E_i}$. Hence, the eigenvalues for the action of left multiplication by $c(S)$ on $R$ are  $\{f_{S,i}\}$, each of them with multiplicity $d_{i}^{2}$, where $d_i:=dim(E_i)$.
\\
\begin{lemma}
    For each irreducible character $\chi_i$ and corresponding codegree $f_i$ we have $f_{S,i}= \chi_i (S) f_i$ for every simple $S$, and so  $\chi_i(c(S))=d_i f_i  \chi_i (S)$.
\end{lemma}
\begin{proof}

For each $\rho_i$, denote by $\chi_i$ the corresponding irreducible character, then we have $\chi_i (c(S))=d_i f_{S,i}$ for each simple $S$. Therefore,
\begin{equation}
	d_i f_{S,i}=\chi_i (c(S))=\operatorname{Tr}(\rho_i (c(S)))=\operatorname{Tr}(\rho_i (\sum_X XSX^{*}))=\operatorname{Tr}((\sum_X \rho_i (XSX^{*}))=
\end{equation}
\begin{equation*}
	=\sum_X \operatorname{Tr}(\rho_i (XSX^{*}))=\sum_X \operatorname{Tr}(\rho_i (X)\rho_i(S)\rho_i(X^{*}))=\sum_X \operatorname{Tr}(\rho_i (S)\rho_i(X^{*})\rho_i(X))=
\end{equation*}

\begin{equation*}
	=\operatorname{Tr}(\rho_i (S)\rho_i(c(1)))=\operatorname{Tr}(\rho_i (S)d_i f_i id_{E_i})=d_i f_i \operatorname{Tr}(\rho_i (S))= d_i f_i  \chi_i (S)
\end{equation*}
\medskip

 because $\operatorname{Tr}(ABC)=\operatorname{Tr}(CAB)=\operatorname{Tr}(BCA)$ holds if $A,B,C$ are square matrices of the same size.
 \end{proof}
 
 \medskip

\begin{corollary}
 	The eigenvalues of the action of $c(S)$ on $R$ are $\{\chi_i (S) f_i\}$, each of them with multiplicity $d_{i}^{2}$. Hence, $det_S=det_1 \prod\limits_{i=1}^{r}
	\left(
	\frac{\chi_i(S)}{d_i}
	\right)^{d_i^2}
	$
\end{corollary}
\bigskip

\section{Character values in terms of the $S$-matrix}
 From the previous section we have $\chi_i(c(S))=d_i f_i  \chi_i (S)$ for each irreducible character $\chi_i$  of $R$ and each simple object $S$ in $\mathcal{C}$. The forgetful functor $F:\mathcal{Z}(\mathcal{C})\rightarrow \mathcal{C}$ induces a ring homomorphism $Gr(\mathcal{Z}(\mathcal{C}))\otimes \mathbb{C}\rightarrow R$ which with an abuse of notation we will also denote by $F$. Then, for each $i$, $\widehat{\rho_i}:=\rho_i \circ F$ is a representation of $Gr(\mathcal{Z}(\mathcal{C}))$, and any object $A$ in $\mathcal{Z}(\mathcal{C})$ acts as $F(A)$ on each $E_i$, i.e., $\widehat{\rho_i}(A)=\rho_i (F(A))$. Furthermore, the image of $F$ is in the center of $R$ (see \cite[Section 7.13]{EGNO}), and so by Schur's lemma we have $\rho_i (F(A))=\alpha_{A,i} ~ id_{E_i}$ for some scalar $\alpha_{A,i}$. Then the assignment $A\rightarrow \alpha_{A,i}$ clearly induces  an irreducible character on $Gr(\mathcal{Z}(\mathcal{C}))$, which we will denote by $\widehat{\chi_i}$. Then we have:
 \begin{equation}
 	\chi_i (F(A))=\operatorname{Tr}(\rho_i (F(A)))=d_i \alpha_{A,i}= d_i \widehat{\chi_i}(A)
 \end{equation}
\medskip

%The induction functor
%$I:\mathcal{C}\to\mathcal{Z}(\mathcal{C})$
%is both a left and a right adjoint to $F$ (cite).

%More precisely, for every %$X\in\mathcal{C}$ and
%$A\in\mathcal{Z}(\mathcal{C})$, there are natural
%isomorphisms\[
%\operatorname{Hom}_{\mathcal{Z}(\mathcal{C})}(I(X),A)
%\cong\operatorname{Hom}_{\mathcal{C}}(X,F(A))
%\]and\[\operatorname{Hom}_{\mathcal{Z}(\mathcal{C})}(A,I(X))
%\cong\operatorname{Hom}_{\mathcal{C}(F(A),X).\]
%The forgetful functor $F$ has a right and left adjoint functor

For each simple object $S$ in $\mathcal{C}$ we have $c(S)=F(I(S))$. In particular, for $A=I(S)$, we get:

\begin{equation}
	\widehat{\chi_i}(I(S))=\frac{\chi_i (c(S))}{d_i}=\frac{d_i f_i \chi_i (S)}{d_i}=f_i \chi_i (S)
\end{equation}

which follows from Lemma 3.2 and $(2)$.

\bigskip

    It is well-known that every irreducible character over  $Gr(\mathcal{Z}(\mathcal{C}))\otimes \mathbb{C}$ has the form 
    
    $$
    \widehat{\chi}_B: A\rightarrow \frac{S_{B,A}}{dim(B)}
    $$ 
     for some simple object $B$ in the center; where if $A$ is simple, $S_{B,A}$ denotes the corresponding entry of the $S$-matrix. We will use the same notation even if $A$ is not simple, meaning, if $A=\sum_j m_j A_j$, where the $A_j$'s are non-isomorphic simples, then $S_{B,A}:=\sum_j m_j S_{B,A_j}$
\medskip

It turns out we can write each character value in terms of the (normalized) S-matrix as follows:

\begin{proposition}
	 Let $\chi_i$ be a irreducible character over $R$. Then for each simple object $S$ in $\mathcal{C}$ we have $\chi_i (S)=\frac{S_{B_i,I(S)}}{dim(\mathcal{C})}$ for some direct summand $B_i$ of $I(1)$. In particular, $S_{B_i,I(1)}=dim(\mathcal{C}) d_i=dim(\mathcal{C})[I(1):B_i]$.
\end{proposition}
\begin{proof}
	Let $\chi_i$ be a irreducible character over $R$ and let $f_i$ be the corresponding codegree. Then we have $\widehat{\chi_i}=\widehat{\chi}_{B_i}$ for some simple direct summand $B_i$ of $I(1)$ (see \cite{O2}). For $A=I(S)$ it follows from $(3)$ that
\begin{equation}
	f_i \chi_i (S)=\widehat{\chi_i}(I(S))=\widehat{\chi}_{B_i}(I(S))=\frac{S_{B_i,I(S)}}{dim(B_i)}
\end{equation}
Furthermore, it follows from \cite[Theorem 2.13]{O2} that $dim (B_i)=\frac{dim(\mathcal{C})}{f_i}$, and so the first part holds.

\medskip 
For the last part, notice that $d_i = \chi_i (1)$, and that $d_i= [I(1):B_i]$ by \cite[Theorem 2.13]{O2}.
\end{proof}

\begin{remark}
As a consequence of the preceding lemma, the irreducible
characters of $K(\mathcal{C})$ take values in a cyclotomic ring of
integers, recovering a standard result in fusion categories $($see for example \cite[Corollary 9.3.13]{EGNO}$)$.
\end{remark}
%\begin{proof}
	%Both $dim(\mathcal{C})$ and any entry of an $S$-matrix are known to be cyclotomic integers(cite). Notice that each $\chi_i (S)$ is an algebraic integer: $\rho_i$ is a direct summand of the regular representation, and so each eigenvalue of $\rho_i (S)$ is an eigenvalue of the fusion matrix $N_S$, which has integer entries. Hence, each eigenvalue must be an algebraic integer, and so must be any finite sum of those. In particular so is $\chi_i (S)=tr(\rho_i (S))$. But then, denoting by $\mathcal{D}$ the set of algebraic integers, and by $N$ the conductor of the center, it follows from the above proposition that $\chi_i (S)\in \mathcal{D}\cap \mathbb{Q}(\zeta_N)=\mathbb{Z}[\zeta_N]$.
%\end{proof}

\medskip

The following lemma shows that the family of induced objects \(\{I(S)\}_{S\in\operatorname{Irr}(\mathcal C)}\) separates the irreducible characters of \(R\).
\begin{lemma}
 Let $A, A'$ be two simple summands of $I(1)$ with corresponding irreducible representations $E,E'$ over $R$. Assume $\widehat{\chi}_A (I(S))=\widehat{\chi}_{A'} (I(S))$ for all simple objects $S$ in $\mathcal{C}$. Then  $A=A'$.
\end{lemma}
\begin{proof}
        Under these conditions we have $f_E \chi_E (S)=f_{E'} \chi_{E'} (S)$ for every simple object. Assume $E\neq E'$. The irreducible characters are orthogonal \cite[Section 8.14]{EGNO}, i.e. $\sum_S \chi_E (S)\chi_{E'} (S^*)=0$. Therefore, $\sum_S f_E\chi_E (S)\chi_{E'} (S^*)=0$, which substituting the expression above gives

        \begin{equation}
            0=\sum_S f_{E'}\chi_{E'} (S)\chi_{E'} (S^*)=f_{E'}^2 dim(E')
        \end{equation}
        which is a contradiction, as every codegree is non-zero.
    \end{proof}
\bigskip

\section{$I(1)$ is Galois closed}

Let $\widehat{\chi}_A$ denote any irreducible character in the center. We know that, if $A$ is a simple direct summand of $I(1)$, then $A$ correspond to an irreducible representation $E$ of $R$, and in that case we have shown in the previous section that $\widehat{\chi}_A (I(S))=f_E \chi_E (S)$ for any simple object $S$ in $\mathcal{C}$. In particular, we have $0\neq dim(E)f_E=\widehat{\chi}_A (I(1))$.
On the other hand, we know that $\widehat{\chi}_A (B)=\frac{S_{A,B}}{dim(A)}$ for any object $B$ in the center. Consequently, by Proposition 4.1, we have $S_{A,I(S)}=D \chi_E (S)$, where $D$ denotes the dimension of $\mathcal{C}$.

%\begin{remark}
 %   This lemma means that 
%\end{remark}

\begin{proposition}
    Let $A$ be a simple representation in the center that is not a direct summand of $I(1)$ then $0=\widehat{\chi}_A (I(S))$, for any simple object $S$ in $\mathcal{C}$. In particular we have $0=\widehat{\chi}_A (I(1))$.
\end{proposition}
\begin{proof}
    Denote the class of $I(1)$ as $[I(1)]$. Then \cite[Proposition 2.8]{O2} tells us that $[I(1)]$ is an eigenvector for the action of the $S$-matrix on $K(\mathcal Z(\mathcal C))\otimes_{\mathbb Z}\mathbb C$, with eigenvalue $D$, where the basis consists of $\{[A] : A \text{ is simple in } \mathcal Z(\mathcal C)\}$. That is, $S\cdot [I(1)]=D[I(1)]$. Hence, we have $S_{A,I(1)}=D[I(1):A]$ for any simple $A$. In particular, $0=S_{A,I(1)}$ for any simple $A$ which is not a direct summand of $I(1)$. And so 
    \begin{equation}
        \widehat{\chi}_A (I(1))=\frac{S_{A,I(1)}}{dim(A)}=0
    \end{equation}
    
For any other simple $S$, notice by \cite[Lemma 8.12.2]{EGNO} that the induced class $[I(S)]$ belongs to the ideal generated by $[I(1)]$ in $\operatorname{Gr}(\mathcal Z(\mathcal C))\otimes\mathbb C$, so the result follows since $\widehat{\chi}_A $ is multiplicative.
    
   \end{proof}
   \medskip

    %\textcolor{blue}{From here we can deduce the following result, which seems to be known by some experts on the matter...but write it in the introduction}

     Let $\sigma$ denote any Galois automorphism of $K=\mathbb{Q}(\zeta_N)$, where $N$ is the conductor. It is well known that the Galois group acts on the irreducible characters of the Grothendieck ring of $\mathcal{Z}(\mathcal{C})$ via
    \begin{equation}
        \sigma\left(\frac{S_{A,B}}{\dim(A)}\right)
=
\frac{S_{A',B}}{\dim(A')}
    \end{equation}
     This induces a permutation of the simple objects of $\mathcal{Z}(\mathcal{C})$, and we write $A'=\sigma(A)$. In terms of the irreducible characters this means $\sigma(\widehat{\chi}_A (B))=\widehat{\chi}_{A'} (B)$ for any object $B$ in $\mathcal{Z}(\mathcal{C})$.
\begin{theorem}
    $I(1)$ is Galois closed, i.e. the action of any Galois automorphism $\sigma$ on a direct summand of $I(1)$ gives another direct summand of $I(1)$.
\end{theorem}
\begin{proof}
Let $A$ a direct summand of $I(1)$, with corresponding irreducible representation $E$ of $R$. Then for $B=I(1)$ we have
    \begin{equation}
        dim(E)\sigma(f_E)=\sigma(\widehat{\chi}_A (I(1)))=\widehat{\chi}_{A'} (I(1))
    \end{equation}
    The LHS is never $0$, and the RHS is $0$ if $A'$ is not a direct summand of $I(1)$ by the proposition above. This forces $A'$ to be another direct summand of $I(1)$.
\end{proof}
\begin{remark}
    It had previously been suggested to the author by Prof. Ng that the statement above is true, while the proof given here was developed independently. A related Galois-invariance result was proved in \cite{K} for Lagrangian algebras in unitary TQFTs. The argument used there can be adapted to the spherical setting, providing an alternative proof of our statement.
\end{remark}
%\begin{remark}
    %A related Galois invariance property for Lagrangian algebras in the unitary setting was observed in [reference]. Our result applies to the Drinfeld center of an arbitrary spherical fusion category.
%\end{remark}

\section{odd near-group categories}

Let \(\mathcal{C}\) be a near-group fusion category of type \((G,m)\), where \(G\) is a finite abelian group and \(m\) is a non-negative integer, with both $
n=|G|$ and $m$ odd. The simple objects of \(\mathcal{C}\) are \(G\cup\{X\}\), with fusion rules
\[
gh=gh,\qquad gX=Xg=X,
\]

and
\[
X^2=\sum_{g\in G}g+mX,
\]

for all \(g,h\in G\).
\medskip

 In this case, the Casimir element is $c(1)=\sum_{S\in\operatorname{Irr}(\mathcal{C})}S^*S
=n1+\sum_{g\in G}g+mX$, and its determinant is $\det_1=n^n(m^2+4n)$\cite{WLL}. This follows from the fact that $n$ is a codegree with multiplicity $n-1$, while the remaining two codegrees are $2n+\frac{m}{2}\left(m\pm\sqrt{m^2+4n}\right)$. In particular we have  
\begin{equation}
    D=FPdim(\mathcal{C})=2n+\frac{m}{2}\left(m+\sqrt{m^2+4n}\right)
\end{equation}

\medskip

Notice that $D$ must be is irrational, since no integral near-group category exists with $G$ abelian and both $n=|G|$ and $m$ odd (see \cite[Lemma A.9 and Theorem A.14]{O2}). Then, because of \cite[Proposition 3.3]{Sch} we have that $n$ always divides $m$, so we can write $m=nk$. Also, $\mathcal{C}$ can be assumed to be spherical and pseudounitary \cite{Th}, so we will write $d=FPdim(X)=dim(X)$, and so $D=n+d^2$.
\medskip

%The first thing to notice is the following:

\begin{lemma}
    The conductor $N$ of $\mathcal{C}$ is odd.
\end{lemma}
\begin{proof}
    This follows from the fact that $\det_1=n^n(m^2+4n)$ is odd, and we know the determinant and the conductor share the same prime factors when the Grothendieck ring is commutative by Lemma 3.1.
\end{proof}
\begin{corollary}
    Let $\omega$ be a primitive $N$th-root of unity, and let $\sigma_2$ the field endomorphism of $K=\mathbb{Q}(\omega)$ that sends $\omega$ to $\omega^2$. Then $\sigma_2$ is a Galois automorphism.
\end{corollary}
    \medskip

 Notice that  $d= \frac{m+\sqrt{m^2+4n}}{2}$, so $L=\mathbb{Q}(d)=\mathbb{Q}(\sqrt{m^2+4n})$; and since $m$ and $n$ are odd, we have $m^2+4n\equiv 5$ (mod $8$). its square free part is also congruent to $5$ (mod $8$). Hence, the prime $2$ is inert in $L$ and consequentely $\sigma_2$ must act non-trivially on $L$. This means in particular that $d':=\sigma_2(d)=\frac{m-\sqrt{m^2+4n}}{2}=m-d$.

\begin{lemma}
    Let $B$ be the direct summand of $I(1)$ that corresponds to the codegree $D'=2n+\frac{m}{2}\left(m-\sqrt{m^2+4n}\right)$, i.e. $dim B=\frac{D}{D'}$. Then $\sigma_2 (\mathbf{1})=B$, where $\mathbf{1}$ denotes the unit object in the center.
\end{lemma}
\begin{proof}
    From the previous section we know that $I(1)$ is Galois closed. Hence, calling $A=\sigma_2 (\mathbf{1})$, we have
\begin{equation}
        D'=\sigma_2(D)=\sigma_2(\widehat{\chi}_{\mathbf{1}} (I(1)))=\widehat{\chi}_{A} (I(1))=f_E
    \end{equation}
    where $E$ is the irreducible representation over $R$ corresponding to $A$. Since $D'$ only appears once as a codegree we necessarily have $A=B$.
\end{proof}

\begin{lemma}
    We have $F(B)=1\oplus kX$.
\end{lemma}
\begin{proof}
   It is easy to see that $\frac{D}{D'}=1+kd$, and so $dim(F(B))=dim(B)=\frac{D}{D'}=1+kd$. Write now $F(B)=\sum_g a_g  +rX$, then taking dimensions we have 
 $dim(F(B)=a+rd$, where $a=\sum_g a_g$. Hence, we get $1+kd=a+rd$, and since $d$ is irrational this forces $r=k$ and $a=1$. Finally, $1=[I(1):B]=[1:F(B)]$, so the result follows. 
\end{proof}
\medskip

For every simple object $Y$ of $\mathcal{Z}(\mathcal{C})$ write
\[
F(Y)=\bigoplus_{g\in G}a_{Y,g}g\oplus b_YX,
\qquad a_Y=\sum_{g\in G}a_{Y,g}.
\]
Then $a_Y,b_Y\in\mathbb Z_{\geq0}$ and
$\dim(Y)=a_Y+b_Yd$.

\begin{lemma}
  For every simple $Y$,
\begin{equation}\label{eq:two-S-rows}
S_{B,Y}=(1+kd)a_Y-db_Y,
\qquad S_{\mathbf 1,Y}=a_Y+db_Y.
\end{equation}
In particular,
\[
S_{B,B}=1.
\]

\end{lemma}

\begin{proof}

For every $Y$ we have
\[
\frac{S_{B,Y}}{\dim(B)}
=\sigma_2\!\left(\frac{S_{\mathbf 1,Y}}{\dim(\mathbf 1)}\right)
=\sigma_2(\dim Y)=a_Y+b_Yd'.
\]
Hence, since $(1+kd)d'=d'-m=-d$, we have
\begin{equation}
    S_{B,Y}=dim(B)(a_Y+b_Yd')=(1+kd)(a_Y+b_Yd')=(1+kd)a_Y-db_Y
\end{equation}
 In particular, for $Y=B$, we have 
\begin{equation}
    S_{B,B}=(1+kd)(1+kd')=1+kd-kd=1
\end{equation}
since $d'd=-n$.
\end{proof}

\medskip

It is well-known that the Drinfeld center of a spherical fusion category has trivial multiplicative central charge (see \cite[Section 8.27]{EGNO}), so its positive Gauss sum is $\tau^+(\mathcal Z)=D$.
With the (unnormalized) $S$-matrix, the modular identity becomes
\begin{equation}\label{eq:STS-k1}
STS=D\,T^{-1}ST^{-1}.
\end{equation}

\begin{proposition}
We have $2n=\sum_Ya_Y^2\theta_Y$,~ $m=\sum_Ya_Yb_Y\theta_Y$,~ $kd'=\sum_Yb_Y^2\theta_Y$, where the sums are over all simples objects $Y$ in the center, and $\theta_Y$ denotes the twist of $Y$.
\end{proposition}
\begin{proof}
    Call \[
R_1=\sum_Ya_Y^2\theta_Y,\qquad
P_1=\sum_Ya_Yb_Y\theta_Y,\qquad
Q_1=\sum_Yb_Y^2\theta_Y.
\]

Since $\theta_{\mathbf 1}=\theta_B=1$ \cite[Remark 4.6]{NS}, taking the
$(\mathbf 1,\mathbf 1)$, $(\mathbf 1,B)$, and $(B,B)$ entries
of \eqref{eq:STS-k1}, respectively, and substituting
\eqref{eq:two-S-rows}, gives
\begin{align}
R_1+2dP_1+d^2Q_1&=D,\label{eq:R1}\\
(1+kd)R_1+kd^2P_1-d^2Q_1&=D(1+kd),\label{eq:R2}\\
(1+kd)^2R_1-2(1+kd)dP_1+d^2Q_1&=D.\label{eq:R3}
\end{align}
Subtracting \eqref{eq:R3} from \eqref{eq:R1} gives
\[
(kd)(2+kd)R_1+2d(2+kd)P_1=0.
\]
Hence, we obtain
\[
P_1=\frac{k}{2}R_1.
\]
Subtracting $1+kd$ times \eqref{eq:R1} from \eqref{eq:R2} gives
\[
-d(2+kd)P_1-d^2(2+kd)Q_1=0,
\]
so $P_1=-dQ_1$. Thus \eqref{eq:R1} becomes
$D=R_1+dP_1=R_1 (1+\frac{kd}{2})$. Since $D=n(2+kd)$, it follows that
\begin{equation}\label{eq:three-twist-sums}
R_1=2n,\qquad P_1=kn=m,\qquad
Q_1=-\frac{m}{d}=kd'.
\end{equation}
\end{proof}
\medskip

We are ready now to prove the main result of this article. 
\begin{theorem}\label{prop:k-equals-one}
Let $\mathcal C$ be a  near-group fusion category of type $G+m$, where $G$ is an abelian group of order $n$, and both $m$ and $n$ are odd. Then $m=n$.
\end{theorem}
\begin{proof}
    
The coefficients $a_Y,b_Y$ are integers, while
$\sigma_2(\theta_Y)=\theta_Y^2$ because each twist is an
$N$th root of unity. Applying $\sigma_2$ to the identities
for $P_1$ and $Q_1$ in \eqref{eq:three-twist-sums} yields
\begin{equation}\label{eq:twist-square-sums}
\sum_Ya_Yb_Y\theta_Y^2=m,
\qquad \sum_Yb_Y^2\theta_Y^2=k\sigma_2(d')=kd.
\end{equation}

It follows from \cite[Theorem 4.1]{NS} that $\nu_2(X)=\frac{1}{D}
\sum_{Y\in\operatorname{Irr}(\mathcal Z)}
\dim(Y)\,[F(Y):X]\,\theta_Y^2$.
\medskip

Since $[F(Y):X]=b_Y$ and $\dim(Y)=a_Y+db_Y$,
\eqref{eq:twist-square-sums} gives
\begin{align*}
D\nu_2(X)
=\sum_Y(a_Y+db_Y)b_Y\theta_Y^2
=m+kd^2=k(n+d^2)=
kD.
\end{align*}
Therefore $\nu_2(X)=k$. But $X$ is simple and self-dual, so $\nu_2(X)\in\{-1,1\}$.
Since $k$ is a positive integer, we conclude that $k=1$; and so $m=n$.

\end{proof}
\medskip

\section*{Acknowledgements}

The author thanks Prof. Ng for useful comments and discussions during the preparation of this work.

\section*{Declaration of AI use}
 The main ideas of this paper are due to the author. GPT-6 Sol (OpenAI) was used to develop the proof of Proposition 6.6. The tool was also used to check some of the arguments, search for relevant literature, and provide references. The author assumes full responsibility for the content of this article.


\begin{thebibliography}{99}

\bibitem[BNRW]{BNRW}
P. Bruillard, S.-H. Ng, E. C. Rowell, and Z. Wang,
\emph{Rank-finiteness for modular categories},
J. Amer. Math. Soc. \textbf{29} (2016), 857--881.

%\bibitem[BG]{BG}
%P. Bruillard, S.-H. Ng, E. C. Rowell, and Z. Wang,
%\emph{On classification of modular categories by rank},
%Int. Math. Res. Not. IMRN (2016), no. 24, 7546--7588.

%\bibitem[EGO]{EGO}
%P. Etingof, S. Gelaki, and V. Ostrik,
%\emph{Classification of fusion categories of dimension $pq$},
%Int. Math. Res. Not. (2004), no. 57, 3041--3056.

\bibitem[EGNO]{EGNO}
P. Etingof, S. Gelaki, D. Nikshych, and V. Ostrik,
\emph{Tensor Categories},
Mathematical Surveys and Monographs, vol. 205,
Amer. Math. Soc., Providence, RI, 2015.

\bibitem[KKOS] {K}
J.~Kaidi, Z.~Komargodski, K.~Ohmori, S.~Seifnashri, and S.-H.~Shao,
\emph{Higher central charges and topological boundaries in $2+1$-dimensional TQFTs},
SciPost Phys. \textbf{13} (2022), no.~3, 067.



%\bibitem[I]{I}
%M. Izumi,
%\emph{A Cuntz algebra approach to the classification of near-group categories},
%Proc. Centre Math. Appl. Austral. Nat. Univ.
%\textbf{46} (2017), 222--343.

%\bibitem[L]{L}
%M. Lorenz,
%\emph{Some applications of Frobenius algebras to Hopf algebras},
%in Groups, Algebras and Applications,
%Contemp. Math., vol. 537,
%Amer. Math. Soc., 2011, 269--289.

\bibitem[NS]{NS}
S.-H. Ng and P. Schauenburg,
\emph{Frobenius--Schur indicators and exponents of spherical categories},
Adv. Math. \textbf{211} (2007), 34--71.

\bibitem[O1]{O1}
V. Ostrik,
\emph{On formal codegrees of fusion categories},
Math. Res. Lett. \textbf{16} (2009), 895--901.

\bibitem[O2]{O2}
V. Ostrik,
\emph{Pivotal fusion categories of rank 3},
Mosc. Math. J. \textbf{15} (2015), 373--396.

\bibitem[Sc1]{Scn}
A. Schopieray,
\emph{Norm, trace, and formal codegrees of fusion categories},
J. Algebra \textbf{568} (2021), 362--385.



\bibitem[Sc2]{Sch}
A. Schopieray,
\emph{Categorification of integral group rings extended by one dimension},
J. Lond. Math. Soc. \textbf{108} (2023), 1617--1641.




\bibitem[Si]{Si}
J. Siehler,
\emph{Near-group categories},
Algebr. Geom. Topol. \textbf{3} (2003), 719--775.

\bibitem[Th]{Th} J. Thornton, \textit{On braided near-group categories}, arXiv:1102.4640 (2011).

\bibitem[TY]{TY}
D. Tambara and S. Yamagami,
\emph{Tensor categories with fusion rules of self-duality for finite abelian groups},
J. Algebra \textbf{209} (1998), 692--707.

\bibitem[WLL]{WLL}
Z. Wang, G. Liu, and L. Li,
\emph{The Casimir number and the determinant of a fusion category},
Glasg. Math. J. \textbf{63} (2021), 438--450.

\end{thebibliography}
\end{document}